\documentclass{amsart}
\usepackage{mathrsfs}
\usepackage{graphicx}
\newtheorem{thm}{Theorem}[section]
\newtheorem{cor}[thm]{Corollary}
\newtheorem{lem}[thm]{Lemma}

\theoremstyle{definition}

\newtheorem{proper}[thm]{Property}
\theoremstyle{remark}
\newtheorem{rem}[thm]{Remark}
\numberwithin{equation}{section}

\begin{document}
\hfill{\textit{Dedicated to Vladimir Mikhailovich Zolotarev,}}

\hfill{\textit{Victor Makarovich Kruglov,}}

\hfill{\textit{ and to the memory of Vladimir Vyacheslavovich
Kalashnikov}}

\title[Continuity theorems for the $M/M/1/n$
queueing system]{Continuity theorems for the $M/M/1/n$
queueing system}%
\author{Vyacheslav M. Abramov}%
\address{School of Mathematical Sciences, Monash University,
Building 28M, Clayton Campus, Clayton, VIC 3800}%
\email{Vyacheslav.Abramov@sci.monash.edu.au}%


\begin{abstract}
In this paper continuity theorems are established for the number
of losses during a busy period of the $M/M/1/n$ queue. We consider
an $M/GI/1/n$ queueing system where the service time probability
distribution, slightly different in a certain sense from the
exponential distribution, is approximated by that exponential
distribution. Continuity theorems are obtained in the form of one
or two-sided stochastic inequalities. The paper shows how the
bounds of these inequalities are changed if further assumptions,
associated with specific properties of the service time
distribution (precisely described in the paper), are made.
Specifically, some parametric families of service time
distributions are discussed, and the paper establishes uniform
estimates (given for all possible values of the parameter) and
local estimates (where the parameter is fixed and takes only the
given value). The analysis of the paper is based on the level
crossing approach and some characterization properties of the
exponential distribution.

\end{abstract}
\maketitle
\section{Introduction}

\subsection{Motivation, the problem formulation and review of the related literature}

In a large number of engineering applications of queueing theory
it is pertinent to know what to expect when we replace one
probability distribution function of a given stochastic model, by
another probability distribution with simpler/concise properties.
Many of these applications are associated with the case where a
probability distribution is replaced by an exponential
distribution if it is known that the given distribution is close
(in a certain sense) to the exponential distribution. These
problems are closely related to characterization problems
associated with the exponential distribution and its stability
(e.g. Azlarov and Volodin \cite{Azlarov and Volodin 1986}).

Continuity analysis of queueing systems is a difficult problem
when compared to general characterization and continuity problems
for random variables, and it is relevant to know the expected
behavior of the system under a range of {\it different}
conditions. The paper discusses three different conditions under
which the probability distribution function of a service time,
which differs slightly from the exponential distribution in the
uniform metric, can be approximated by that exponential
distribution.

This paper considers an $M/GI/1/n$ queueing system, where $n$ is
the buffer size excluding a customer in service, $\lambda$ is the
rate of Poisson input, $B(x)=\mathrm{Pr}\{\chi\le x\}$ is the
probability distribution of a service time $\chi$, and the
parameter $\mu$ is the reciprocal of the expected service time. It
is assumed then that the probability distribution function $B(x)$
is close (under specific cases described below) to the exponential
distribution with the same parameter $\mu$.

Let $L_n$ denote the number of losses during a busy period of this
queueing system. The classic explicit results for the $M/GI/1/n$
queueing system, for example the recurrence relations for the
expectations $\mathrm{E}L_n$, as well as a number of results on
the asymptotic behavior of $\mathrm{E}L_n$ as $n\to\infty$, can be
found in Abramov \cite{Abramov - book 1991}, \cite{Abramov 1997},
Cooper and Tilt \cite{Cooper and Tilt 1976} and Tomko \cite{Tomko
1967}. An application of these explicit results to asymptotic
analysis of lost messages in communication networks is given in
Abramov \cite{Abramov 2004}, and their application to problems of
optimal control of large dams is given in \cite{Abramov 2007} and
\cite{Abramov 2007a}. Some simple stochastic inequalities for
different finite buffer queueing systems are obtained in Abramov
\cite{Abramov 2001a}, \cite{Abramov 2005}, Righter \cite{Righter
1999}, Pek\"oz {\it et al.} \cite{Pekoz et al 2003}, Wolff
\cite{Wolff 2002} and other papers. Other relevant studies of the
loss probability in finite buffer systems can be found in Abramov
\cite{Abramov 1991}, \cite{Abramov 2002} and Choi {\it et al.}
\cite{Choi et al 2000}.

The present paper establishes continuity theorems in the form of
\textit{stochastic inequalities} for $L_n$ similar to the
stochastic inequalities established in Abramov \cite{Abramov
 2001a}, \cite{Abramov 2005}. The main difference between the results of this
paper and those of Abramov \cite{Abramov
 2001a}, \cite{Abramov 2005} is that
the stochastic inequalities established in the present paper
depend on specific properties of the service time distribution,
whereas the stochastic bounds of the earlier papers \cite{Abramov
2001a}, \cite{Abramov 2005} are obtained without any special
assumption on that service time distribution. Then the specific
feature of the present paper is that the stochastic inequalities
obtained can be considered as \textit{approximation bounds} or
\textit{continuity bounds} for the queueing system in the uniform
metric, and these approximation bounds can have real application
for stochastic systems arising in practice. The phrase
\textit{continuity theorems} used in this paper is associated with
the traditional ````$\epsilon$-$\delta$ language" of mathematical
analysis". The input characteristic of the system is assumed to
have a small variation $\epsilon$, and then the output
characteristic has a small variation $\delta$ continuously
depending on $\epsilon$. Both these $\epsilon$ and $\delta$
variations are assumed to be described by some special metrics. A
large variety of metrics for the continuity analysis of complex
systems can be found in Dudley \cite{Dudley 1976}, Kalashnikov and
Rachev \cite{Kalashnikov and Rachev 1990}, Rachev \cite{Rachev
1991}, Zolotarev \cite{Zolotarev 1983} and many others. In the
cases where an application of the well-known traditional metrics
such as the uniform metric or the Levy metric becomes unavailable,
very difficult or not profitable, the continuity analysis of a
complex system requires to use the special type of metrics,
appropriate for analysis of a given system or associated with
given special conditions.

The continuity analysis of queueing systems goes back to the
papers of Kennedy \cite{Kennedy 1972} and Whitt \cite{Whitt 1974}
and further developed by Kalashnikov \cite{Kalashnikov 1983},
Zolotarev \cite{Zolotarev 1976}, \cite{Zolotarev 1977} and
Gordienko and Ruiz de Ch\'avez \cite{Gordienko and Ruiz de Chavez
1998}, \cite{Gordienko and Ruiz de Chavez 2005} to mention a few.

\smallskip

The continuity analysis of such standard characteristics as
\textit{the number of losses during a busy period} is a very
difficult problem. One of the advantages of the method in this
paper, is that the continuity of this characteristic is studied in
terms of the \textit{uniform probability metric}. As a result, the
final results are simple and clear, and the conditions of the
theorems are verifiable. (The uniform metric of one-dimensional
distributions is usually called the Kolmogorov metric.)
Furthermore, under different assumptions that the service time
distribution belongs to its parametric family of distributions
depending on parameter $p$ we establish {\it uniform} estimates
for the number of losses during a busy period given for all
possible values $p$. That is, stochastic inequalities depending
only on $\epsilon$ but not on $p$ are obtained. In the case where
parameter $p$ of this family is fixed (specifically the case
$p=\frac{1}{2}$ is considered) we establish essentially stronger
estimates than these corresponding uniform estimates.

For the purpose of the continuity analysis, this paper essentially
develops the level crossing approach of the earlier papers of
Abramov \cite{Abramov - book 1991}, \cite{Abramov 1994},
\cite{Abramov 2001a}, \cite{Abramov 2001b}, \cite{Abramov 2005}
and some characterization theorems associated with the exponential
distribution, that are then used together with the aforementioned
metrical approach. We also develop the continuity theorem of
Azlarov and Volodin \cite{Azlarov and Volodin 1986} under some
additional conditions and use the properties of the known classes
of special probability distributions such as NBU (New Better than
Used), NWU (New Worse than Used) (e.g. Stoyan \cite{Stoyan 1983}
for the properties of these classes of distribution and related
results). The special parametric order relation
$\mathscr{C}_{\lambda}$, which was introduced and originally
studied in Abramov \cite{Abramov - book 1991}, is also used. For
the reader's convenience all necessary facts and concepts are
recalled in this paper.

\subsection{The main conditions for the probability distribution of
a service time} The main conditions for the probability
distribution function $B(x)$ ($B(0)=0$) of a service time under
which we obtain the stochastic bounds are as follows.

\medskip
$\bullet$ \textit{Condition (A)}. The probability distribution of
a service time has the representation
\begin{equation}\label{r10}
B(x)=pF(x)+(1-p)(1-\mbox{e}^{-\mu x}), \ 0<p\leq 1,
\end{equation}
where $F(x)=\mathrm{Pr}\{\zeta\leq x\}$ is a probability
distribution function of a nonnegative random variable having the
expectation $\frac{1}{\mu}$, and
\begin{equation}\label{r20}
\sup_{x,y\geq 0}\big|F_y(x)-F(x)\big|\leq\epsilon, \ \epsilon\geq
0,
\end{equation}
where $F_y(x)=\mathrm{Pr}\{\zeta\leq x+y~|~\zeta>y\}$. Relation
(\ref{r20}) says that the distance in the uniform metric between
$F(x)$ and $E_\mu(x)=1-\mbox{e}^{-\mu x}$ is not greater than
$\epsilon$, the case $\epsilon=0$ corresponds to the equality
$F(x)=E_\mu(x)$ for all $x$. In the following we assume that
$\epsilon>0$ and write the strong inequalities (i.e. the
right-hand side of (\ref{r20}) is less than $\epsilon$).

\medskip
$\bullet$ \textit{Condition (B)}. Along with (\ref{r10}) and
(\ref{r20}) it is assumed that $F(x)$ belongs either to the class
NBU or to the class NWU.

\smallskip
Recall that a probability distribution function $\Xi(x)$ of a
nonnegative random variable is said to belong to the class NBU if
for all $x\ge 0$, $y\ge 0$ we have $\overline{\Xi}(x+y)\leq
\overline{\Xi}(x)\overline{\Xi}(y)$. (Throughout the paper, for
any probability distribution function $\Xi(x)$ we use the notation
$\overline{\Xi}(x)=1-\Xi(x)$.) If the opposite inequality holds
i.e. $\overline{\Xi}(x+y)\geq \overline{\Xi}(x) \overline{\Xi}(y)$
then $\Xi(x)$ is said to belong to the class NWU.

\medskip
$\bullet$ \textit{Condition (C)}. The probability distribution
function $B(x)$ of a service time belongs to the class NBU, and
\begin{equation}\label{r30}
\sup_{x,y\ge 0}\big|B_y(x)-B(x)\big|<\epsilon,
\end{equation}
where $B_y(x)=\mathrm{Pr}\{\chi\leq x+y~|~\chi>y\}$.

\medskip
$\bullet$ \textit{Condition (D)}. Let
$F_y(x)=\mathrm{Pr}\{\zeta\le x+y~|~\zeta>y\}$ be a family of
given probability distributions. It is assumed that there exists
the probability distribution $F_{y^0}(x)$ of this family
satisfying the relation $F_{y^0}\leq_{\mathscr{C}_{\lambda}}F_y$
for all $y\ge 0$.

The definition and main property of the parametric order relation
$\mathscr{C}_{\lambda}$ is recalled in Appendix A. For a more
detailed consideration see Abramov \cite{Abramov - book 1991}.

\medskip
$\bullet$ \textit{Condition (E)}.
 Let $B_y(x)=\mathrm{Pr}\{\chi\le x+y~|~\chi>y\}$ be a family of given
probability distributions. It is assumed that there exists the
probability distribution $B_{y^0}(x)$ of this family satisfying
the relation $B_{y^0}\leq_{\mathscr{C}_{\lambda}}B_y$ for all
$y\ge 0$.

If the family of probability distributions $F_y(x)$ (or $B_y(x)$)
is partially ordered with respect to the order relation
$\mathscr{C}_{\lambda}$, then according to Zorn's lemma (e.g.
Ciesielski \cite{Ciesielski 1997}) $F_{y^0}$ (or correspondingly
$B_{y^0}$) is a minimal element of this family (regarding this
order relation $\mathscr{C}_{\lambda}$).
\smallskip

The only difference between Conditions (D) and (E) is that
Condition (E) is related to the probability distribution function
$B(x)$ of a service time, whereas Condition (D) is related to the
associated probability distribution function $F(x)$.

Conditions (A), (B) and (C) are the main conditions for our
consideration, while Conditions (D) and (E) are additional
(associated) conditions. In other words one of Conditions (A), (B)
and (C) are always present in our consideration. Condition (D) can
be present only with Conditions (A) and (B), and Condition (E) can
be present only with Condition (C).

It is worth noting that in certain known cases the minimal element
$F_{y^0}$ can be determined easily. For example, in the case where
$F(x)$ belongs to the class NWU, the minimal element of the family
of probability distributions $F_y(x)$ is $F_{y^0}=F_{0}=F$ in the
sense of the order relation $\mathscr{C}_\lambda$. (By replacing
NWU with NBU and respectively `minimal' with `maximal', this fact
is explicitly used in Theorem \ref{t42}, Section 4.)

\subsection{Further discussion of the main conditions}
Conditions (A), (B) and (C) are three different conditions where
if $\epsilon$ is small, then the probability distribution function
$B(x)$ is close to the exponential distribution in the sense of
the uniform metric. Therefore it is very significant to know what
one can expect if Conditions (A), (B) or (C) are satisfied.

In most queueing problems representation (\ref{r10}) is not
standard. The standard assumption appearing in characterization
problems associated with the exponential distribution is
(\ref{r30}) (e.g. Azlarov and Volodin \cite{Azlarov and Volodin
1986}; cf. Daley \cite{Daley 1976} for some earlier results).
Representation of the probability distribution $B(x)$ in its
special form (\ref{r10}) matches the usual case $B(x)=F(x)$ (that
is the case $p=1$). It is shown below (see Lemma \ref{l21}) that
assumptions (\ref{r10}) and (\ref{r20}) of Condition (A) allow us
only to prove  that for all $0<p\leq 1$
\begin{equation}\label{r1_1}
\sup_{x,y\ge 0}\big|B_y(x)-B(x)\big|<5\epsilon,
\end{equation}
and only in the special case where $p=\frac{1}{2}$ we have
\begin{equation}\label{r1_2}
\sup_{x,y\ge 0}\big|B_y(x)-B(x)\big|<2\epsilon
\end{equation}
(see Remark \ref{rem1}). As a result, the approximation bounds
under Conditions (A), (B) and (C) are all different, and more
specifically representation (\ref{r10}) together with assumption
(\ref{r20}) lead to relatively worse estimates than those in the
special case $B(x)=F(x)$ under assumption (\ref{r30}). However
(\ref{r1_1}) is a uniform estimate for all $0<p\leq 1$, whereas
(\ref{r1_2}) is a usual local estimate better than (\ref{r1_1}).
In the special case $p=1$ it is a local estimate too, which in
fact is (\ref{r30}).

The advantage of (\ref{r10}) is also as follows. When the
parameter $p$ is small, the service time distribution is close to
the exponential one in the uniform metric. Thus, our results are
associated with a wider class of approximations and based on a
two-parameter family of distributions. In practice, this can help
us find a more appropriate representation and approximation for
the initial probability distribution function $B(x)$ given
empirically. The main results of this paper establish explicit
dependence of the parameter $\epsilon$ only and are uniform in
$p$. However there are examples related to the special case of
$p=\frac{1}{2}$ where the bounds obtained are essentially better
than in the general case of arbitrary $p$. In a similar manner the
case of any given $p$ can be considered, where the specific
estimates are expected to be better than the corresponding uniform
estimates.

Notice also that the class of stochastic inequalities is wider. On
the one hand, only the parameter $\epsilon$ can be assumed to be
small, and on the other we can assume that both $\epsilon$ and $p$
are small values. The variety of these assumptions enables us to
choose the possible appropriate parameters of the model for the
further approximations in order to obtain a relevant conclusion.

\subsection{Organization and methodology of this paper}
This paper is structured into 4 sections. Following this
introduction, Sections 2, 3 and 4 study the number of losses
during a busy period under Conditions (A), (B) and (C),
respectively.

In Section 2 we give a deeper analysis of the intervals obtained
under the special procedure of deleting subintervals and
connecting the ends as is explained in a number of the earlier
papers of the author (see Abramov \cite{Abramov - book 1991},
\cite{Abramov 1994}, \cite{Abramov 2001a}, \cite{Abramov 2001b},
\cite{Abramov 2005}). Specifically we prove the following results.
The first result of this section, Lemma \ref{l21}, provides the
uniform estimate (\ref{r1_1}) for the parametric family of
distributions. Remark \ref{rem1} provides local estimate
(\ref{r1_2}) in the special case $p=\frac{1}{2}$. The proof of
both these estimates is based on the characterization theorem of
Azlarov and Volodin \cite{Azlarov and Volodin 1986}, which is
formulated in Lemma \ref{l22}. The proof of Theorem \ref{t24} is
based on the level crossing approach which in its present form is
originated and developed by the author (see Abramov \cite{Abramov
- book 1991}, \cite{Abramov 1994}, \cite{Abramov 2001a},
\cite{Abramov 2001b}, \cite{Abramov 2005}).

In Section 3, under Conditions (B) and (D) stronger results than
those under Conditions (A) and (D) are established. The results
are based on similar proofs to those in Section 2.

In Section 4, under Condition (C) the two-sides stochastic
inequalities are given in Theorem \ref{t42}. The proofs of the
statements under Condition (C) are completely analogous to the
earlier proofs under Condition (A).

The partial order relation $\mathscr{C}_\lambda$ and its main
property are discussed in Appendix A. The proofs of Lemma
\ref{l22} and Lemma \ref{l31} is given in Appendix B. In Appendix
C, the proof of Lemma \ref{tex-lem} is given.

\section{The number of losses under Condition (A)}

Throughout the paper we use Kolmogorov's (uniform) metric between
two one-dimensional probability distributions. Recall the
definition of Kolmogorov's metric (e.g. Kalashnikov and Rachev
\cite{Kalashnikov and Rachev 1990}, Rachev \cite{Rachev 1991}).
Let $G(x)=\mathrm{Pr}\{\xi\leq x\}$ and
$H(x)=\mathrm{Pr}\{\eta\leq x\}$ be probability distribution
functions of the random variables $\xi$ and $\eta$. Kolmogorov's
metric $\mathscr{K}(G, H)$ between two probability distribution
functions $G(x)$ and $H(x)$ is defined by
$$
\mathscr{K}(G,H)=\sup_{x\in \mathbb{R}^1}\big|G(x)-H(x)\big|.
$$
In the following for Kolmogorov's metric the notation
$\mathscr{K}(\xi,\eta)$ is also used, where $\mathscr{
K}(\xi,\eta)=\mathscr{K}(G,H)$. All exponential distributions are
denoted $E_\alpha=E_\alpha(x)=1-\mbox{e}^{-\alpha x}$, where
$\alpha>0$ is the parameter of distribution.

\begin{lem}\label{l21}
Under assumptions (\ref{r10}) and (\ref{r20}) we have
$$
\sup_{x\ge 0, y\ge 0}\big|B_y(x)-B(x)\big|< 5\epsilon.
$$
\end{lem}

In order to prove this lemma we use the following result of
Azlarov and Volodin \cite{Azlarov and Volodin 1986}, the proof of
which is also provided in Appendix B.

\begin{lem}\label{l22}
(Azlarov and Volodin \cite{Azlarov and Volodin 1986}.) Let
$F(x)=\mathrm{Pr}\{\zeta\le x\}$ be a probability distribution
function of a nonnegative random variable with the expectation
$1/\mu$, and
$$
\sup_{x\ge 0, y\ge 0}\big|F_y(x)-F(x)\big|<\epsilon.
$$
Then
$$
\mathscr{K}(F,E_\mu)< 2\epsilon.
$$
\end{lem}

\textit{Proof of Lemma 2.1.} We have
\begin{equation}\label{r2.4-1}
\begin{aligned}
\overline{B}_y(x)
&=\frac{1-pF(x+y)-(1-p)E_{\mu}(x+y)}{1-pF(y)-(1-p)E_{\mu}(y)}\\
&= \frac{ p\overline{F}(x+y)+(1-p)\mbox{e}^{-\mu(x+y)}}{
p\overline{F}(y)+(1-p)\mbox{e}^{-\mu y}}.
\end{aligned}
\end{equation}
Therefore,
\begin{equation}\label{r2.4}
\begin{aligned}
&|B_y(x)-B(x)|\\
&=\Big|\frac{ p\overline{F}(x+y)+(1-p)\mbox{e}^{-\mu(x+y)}}{
p\overline{F}(y)+(1-p)\mbox{e}^{-\mu
y}}-  p\overline{F}(x)-(1-p)\mbox{e}^{-\mu x}\Big|\\
&=\Big|\frac{p\overline{F}(x+y)+(1-p)\mbox{e}^{-\mu(x+y)}}{
p\overline{F}(y)+(1-p)\mbox{e}^{-\mu
y}} -\frac{p^2 \overline{F}(x) \overline{F}(y)}{ p\overline{F}(y)+(1-p)\mbox{e}^{-\mu y}}\\
&\ \ \ -\frac{ p(1-p)[\mbox{e}^{-\mu
x}\overline{F}(y)+\mbox{e}^{-\mu y} \overline{F}(x)]} {
p\overline{F}(y)+(1-p)\mbox{e}^{-\mu y}} -\frac{
(1-p)^2\mbox{e}^{-\mu (x+y)}}{
p\overline{F}(y)+(1-p)\mbox{e}^{-\mu y}}\Big|.
\end{aligned}
\end{equation}

According to (\ref{r20}) we have
\begin{equation}\label{r2.5}
\Big|\frac{ \overline{F}(x+y)}{ \overline{F}(y)}-
\overline{F}(x)\Big|<\epsilon,
\end{equation}
where the denominator in the fraction of (\ref{r2.5}) is assumed
to be positive, and if it is equal to zero, then we use the
convention $\frac{0}{0}=0$. It follows from Lemma \ref{l22} that
\begin{equation}\label{r2.6}
| \overline{F}(x)-\mbox{e}^{-\mu x}|<2\epsilon.
\end{equation}
Therefore from (\ref{r2.4}), (\ref{r2.5}), (\ref{r2.6}) and the
triangle inequality for all $x\ge 0$, $y\ge 0$, we obtain:
\begin{equation}\label{r2.7}
\begin{aligned}
&\big|p \overline{F}(x+y) -p^2 \overline{F}(x) \overline{F}(y)
-p(1-p)\mbox{e}^{-\mu
x} \overline{F}(y) \big|\\
&\leq p \overline{F}(y) \Big|\frac{ \overline{F}(x+y)}{
\overline{F}(y)}-\mbox{e}^{-\mu
x}\Big|+p^2 \overline{F}(y)\Big| \overline{F}(x)-\mbox{e}^{-\mu x}\Big|\\
&= p \overline{F}(y) \Big|\frac{ \overline{F}(x+y)}{
\overline{F}(y)}- \overline{F}(x)+ \overline{F}(x)-\mbox{e}^{-\mu
x}\Big|\\
&\ \ \ +p^2 \overline{F}(y) \big| \overline{F}(x)-\mbox{e}^{-\mu x}\big|\\
&\leq p \overline{F}(y) \Big(\Big|\frac{ \overline{F}(x+y)}{
\overline{F}(y)}- \overline{F}(x)\Big|
+\big| \overline{F}(x)-\mbox{e}^{-\mu x}\big|\Big)\\
&\ \ \ +p^2 \overline{F}(y) \big| \overline{F}(x)-\mbox{e}^{-\mu x}\big|\\
&<3\epsilon p \overline{F}(y) +2\epsilon p^2 \overline{F}(y)\\
&\leq 5\epsilon p \overline{F}(y),
\end{aligned}
\end{equation}
\begin{equation}\label{r2.8}
\begin{aligned}
&\big|(1-p)\mbox{e}^{-\mu(x+y)}-p(1-p)\mbox{e}^{-\mu
y}  \overline{F}(x) -(1-p)^2\mbox{e}^{-\mu(x+y)}\big|\\
 &\leq
(1-p)\mbox{e}^{-\mu y}\Big(\big|\mbox{e}^{-\mu x}-
\overline{F}(x)\big|+\big|(1-p) \overline{F}(x)
-(1-p)\mbox{e}^{-\mu
x}\big|\Big)\\
&< 2\epsilon(1-p)\mbox{e}^{-\mu y}+2\epsilon(1-p)^2\mbox{e}^{-\mu
y}\\
&< 3\epsilon(1-p)\mbox{e}^{-\mu y}+2\epsilon(1-p)^2\mbox{e}^{-\mu
y}\\ &< 5\epsilon(1-p)\mbox{e}^{-\mu y},
\end{aligned}
\end{equation}
and then it follows from (\ref{r2.7}) and (\ref{r2.8}) that for
all $x\ge 0$, $y\ge 0$,
$$
|B_y(x)-B(x)|<\frac{5\epsilon p \overline{F}(y)
+5\epsilon(1-p)\mbox{e}^{-\mu y}}{
p\overline{F}(y)+(1-p)\mbox{e}^{-\mu y}}=5\epsilon.
$$
Hence, Lemma \ref{l21} is proved.

\begin{rem}\label{rem1}
The statement of Lemma \ref{l21} is a uniform estimate for all
$p$. In special cases where $p$ is given, one can obtain stronger
estimates. For example, in the case $p=\frac{1}{2}$ we have the
following inequalities:
\begin{equation}\label{r2.7_1}
\begin{aligned}
&\big|\textstyle{\frac{1}{2}} \overline{F}(x+y)
-\textstyle{\frac{1}{4}} \overline{F}(x) \overline{F}(y)
-\textstyle{\frac{1}{4}}\mbox{e}^{-\mu
x} \overline{F}(y) \big|\\
&\leq
\textstyle{\frac{1}{4}} \overline{F}(y) \big|\frac{ \overline{F}(x+y)}{ \overline{F}(y)}-\overline{F}(x)\big|\\
&\ \ \ +\textstyle{\frac{1}{4}} \overline{F}(y) \big|\frac{
\overline{F}(x+y)}{ \overline{F}(y)}-\mbox{e}^{-\mu
x}\big|\\
&\leq \textstyle{\frac{1}{2}} \overline{F}(y) \big|\frac{
\overline{F}(x+y)}{
\overline{F}(y)}- \overline{F}(x)\big|\\
&\ \ \ +\textstyle{\frac{1}{4}} \overline{F}(y) \big| \overline{F}(x)-\mbox{e}^{-\mu x}\big|\\
&<\textstyle{\frac{1}{2}}\epsilon \overline{F}(y)
+\textstyle{\frac{1}{2}}\epsilon \overline{F}(y)\\ &=\epsilon
\overline{F}(y),
\end{aligned}
\end{equation}
\begin{equation}\label{r2.8_1}
\begin{aligned}
\big|\textstyle{\frac{1}{2}}\mbox{e}^{-\mu(x+y)}-\textstyle{\frac{1}{4}}\mbox{e}^{-\mu
y} \overline{F}(x)
-\textstyle{\frac{1}{4}}\mbox{e}^{-\mu(x+y)}\big|
&=\textstyle{\frac{1}{4}}\mbox{e}^{-\mu y}\big|\mbox{e}^{-\mu
x}- \overline{F}(x)\big|\\
&<\textstyle{\frac{1}{2}}\epsilon\mbox{e}^{-\mu
y}\\&<\epsilon\mbox{e}^{-\mu y}.
\end{aligned}
\end{equation}
Therefore, from (\ref{r2.7_1}) and (\ref{r2.8_1}) we obtain:
\begin{equation}\label{r2.8_1-2}
|B_y(x)-B(x)|<\frac{\epsilon
 \overline{F}(y) +\epsilon\mathrm{e}^{-\mu
y}}{ \frac{1}{2}\overline{F}(y)+\frac{1}{2}\mathrm{e}^{-\mu
y}}=2\epsilon.
\end{equation}
\end{rem}

\begin{thm}\label{t24}
Under Conditions (A) and (D) we have
$$
L_n\geq_{st}\sum_{i=1}^{Z_{n}}\varsigma_i.
$$
$Z_n$ denotes the number of offspring in the $n$th generation of
the Galton--Watson branching process with $Z_0=1$ and the
offspring generating function
\begin{equation}\label{r2.8.1add}
g(z)<\frac {\widehat B(\lambda)}{1-z+z\widehat
B(\lambda)}+5\epsilon\cdot \frac {1-z+2z\widehat
B(\lambda)}{1-z+z\widehat B(\lambda)}, \ \ |z|\leq 1,
\end{equation}
where
\begin{equation*}
\begin{aligned}
&\widehat B(\lambda)=p\widehat
F(\lambda)+(1-p)\frac{\mu}{\mu+\lambda},\\
&\widehat F(\lambda)=\int_0^\infty\mathrm{e}^{-\lambda
x}\mathrm{d}F(x).
\end{aligned}
\end{equation*}
 $\varsigma_i$, $i=1,2,\ldots$, is the sequence of nonnegative
integer-valued independent identically distributed random
variables all having the probability law
\begin{equation}\label{r2.8.1add'}
\mathrm{Pr}\{\varsigma_i=m\}=\int_0^\infty\mathrm{e}^{-\lambda
x}\frac{(\lambda x)^m}{m!}\mathrm{d} G(x),
\end{equation}
and $G(x)$ is a probability distribution that satisfies the
inequalities $G\leq_{\mathscr{C}_\lambda} F_{y^0}$ and
$G\leq_{\mathscr{C}_\lambda} E_\mu$, where
$F_{y^0}\leq_{\mathscr{C}_\lambda}F_{y}$ for all $y\geq0$.
\end{thm}

\begin{rem}\label{rem2.6}
The condition \eqref{r2.8.1add'} of Theorem \ref{t24} involves a
probability distribution function $G(x)$. If Condition (D) is not
satisfied, or the probability distribution $F_{y^0}(x)$ is unknown
and cannot be evaluated, then a probability distribution function
$G(x)$ should be chosen such that $G\leq_{\mathscr{C}_\lambda}
F_y$ for all $y\geq0$ and $G\leq_{\mathscr{C}_\lambda} E_\mu$.
\end{rem}

\medskip
\begin{proof} Consider a busy period of an $M/GI/1/n$ queueing system.
Let $f(j), 0\le j\le n+1$, denote the number of customers arriving
during a busy period who, upon their arrival, meet $j$ customers
in the system. It is clear that $f(0)=1$ with probability 1. Let
$t_{j,1}, t_{j,2},...$, $t_{j,f(j)}$ be the instants of arrival of
these $f(j)$ customers, and  let $s_{j,1}, s_{j,2},...$, $s_{j,
f(j)}$ be the instants of service completions (the case $j\leq n$)
or losses (the case $j=n+1$) at which there remain only $j$
customers in the system. Note that $t_{n+1,k}=s_{n+1,k},~~ 1\leq
k\leq f(n+1)$, and $f(n+1)=L_n$, the number of losses during a
busy period.

For $0\le j\le n$ consider the intervals:
\begin{equation}\label{r2.16}
\big[t_{j,1},s_{j,1}\big), ~\big[t_{j,2},s_{j,2}\big),...,
\big[t_{j,f(j)}, s_{j,f(j)}\big).
\end{equation}
 It is clear that the intervals
\begin{equation}\label{r2.17}
\big[t_{j+1,1},s_{j+1,1}\big), ~\big[t_{j+1,2},s_{j+1,2}\big),...,
\big[t_{j+1,f(j+1)}, s_{j+1,f(j+1)}\big)
\end{equation}
are contained in intervals (\ref{r2.16}). Let us delete the
intervals in (\ref{r2.17}) from those in (\ref{r2.16}) and connect
the ends, that is every point $t_{j+1,k}$ with the corresponding
point $s_{j+1,k},~1\le k\le f(j+1)$, if the set of intervals
(\ref{r2.17}) is not empty. In other words, in the interval of the
form $[t_{j,k}, s_{j,k})$, the inserted points are of the form
$t_{j+1,m}$. The random variable $\xi_{j,k}$ is then the length of
interval [$t_{j,k}, s_{j,k}$) minus the intervals $[t_{j+1,m},
s_{j+1,m})$ contained in this interval.

A typical example of the level crossings on a busy period is given
in Figure 1. The arc braces in the figure indicate the places of
connection of the points after deleting the intervals.

\begin{figure}
\includegraphics[width=15cm, height=20cm]{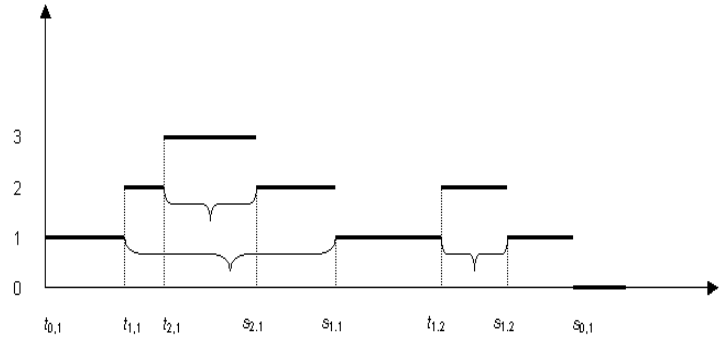}
\caption{Up- and down-crossings on a busy period}
\end{figure}

In the example, given in Figure 1, $n=3$ and the number of
inserted points of level 1 is equal to 2, that is $f(1)=2$.
Similarly, $f(2)=1$.
\medskip

Let us take one of the intervals, [$t_{j,k}, s_{j,k}$) say, and a
customer in service at time $t_{j,k}$. Let $\tau_{j,k}$ denote the
time elapsed from the moment of a service start for that customer
until time $t_{j,k}$, and let $\vartheta_{j,k}$ denote the
residual service time of the tagged customer. The analysis of a
residual service time is well-known in the literature, and in the
simplest cases is associated with the \textit{inspection paradox}
or \textit{waiting time paradox} (see e.g. \cite{Cooper Niu and
Srinivasan 1998}, \cite{Feller 1971}, \cite{Gakis and Sivaslian
1993}, \cite{Herff Jochems Kamps 1997}, \cite{Kremers 1988},
\cite{Resnick 1992}, \cite{Ross 2003} and many others).

The sum $\tau_{j,k}+\vartheta_{j,k}$ has a more complicated
structure than that in the standard inspection or waiting time
paradox, and deriving an explicit representation for
$\mathrm{Pr}\{\vartheta_{j,k}\leq x\}$ in terms of the probability
distributions $B(x)$ or $B_y(x)$ is very difficult. By using Lemma
\ref{l21} above, the probability distribution function
$B_{j,k}(x)=\textrm{Pr}\{\vartheta_{j,k}\leq x\}$ can be evaluated
nevertheless.

Note first the following technical result.
\begin{lem}\label{tex-lem}
For any integer $j\geq1$ and $k\geq1$, and any nonnegative real
$x$ and $y$,
$$
\mathrm{Pr}\{\vartheta_{j,k}\leq x~|~\tau_{j,k}=y\}=B_y(x).
$$
\end{lem}

The proof of this lemma is given in Appendix C.

Therefore, from Lemma \ref{tex-lem} and Lemma \ref{l21} for any
nonnegative $y$ we have:
\begin{equation}\label{r2.18.new}
|\mathrm{Pr}\{\vartheta_{j,k}\leq
x~|~\tau_{j,k}=y\}-B(x)|=|B_y(x)-B(x)|<5\epsilon.
\end{equation}
We will show below, that by using the formula for the total
probability one arrive at
\begin{equation}\label{r2.18.1.new}
|\mathrm{Pr}\{\vartheta_{j,k}\leq x\}-B(x)|<5\epsilon.
\end{equation}
For this purpose we should show that the random variable
$\tau_{j,k}$ is proper, i.e. $\mathrm{Pr}\{\tau_{j,k}<\infty\}=1$,
and then use the formula for the total probability:
\begin{equation}\label{r2.18}
\begin{aligned}
\mathrm{Pr}\{\vartheta_{j,k}\leq x\} &= \int_0^\infty
\mathrm{Pr}\{\vartheta_{j,k}\leq
x~|~\tau_{j,k}=y\}\mbox{d}\mathrm{Pr}\{\tau_{j,k}\leq y~\}\\
&=\int_0^\infty B_y(x)\mbox{d}\mathrm{Pr}\{\tau_{j,k}\leq y~\}
\end{aligned}
\end{equation}
(Notice, that the random variables $\vartheta_{j,k}$ and
$\tau_{j,k}$ are generally dependent, and we cannot derive a
representation for $\mathrm{Pr}\{\tau_{j,k}\leq y\}$.)

Let us show that the random variable $\tau_{j,k}$ is proper for
any $j$ and $k$. Indeed, let $\mathcal{A}_{i}$ denote the event
that there are $i$ arrivals during the time interval from the
service beginning of the customer, who is currently being served
at the moment $t_{j,k}$, until the time moment before the arrival
at $t_{j,k}$ (the arrival at time $t_{j,k}$ and possible arrivals
at the moment of the service began are excluded). Then,
\begin{eqnarray}
\mathrm{Pr}\{\tau_{j,k}\leq
x\}&=&\sum_{i=0}^{j-1}\mathrm{Pr}\{\tau_{j,k}\leq
x~|~\mathcal{A}_{i}\}\mathrm{Pr}\{\mathcal{A}_{i}\}\label{*1}\\
&\geq&\sum_{i=0}^{j-1}\mathrm{Pr}\{\mbox{sum~of} \ i+1 \
\mbox{interarrival~times}\ \leq
x\}\mathrm{Pr}\{\mathcal{A}_{i}\}\label{*2}\\
&\geq&\mathrm{Pr}\{\mbox{sum~of} \ j \ \mbox{interarrival~times}\
\leq x\}.\label{*3}
\end{eqnarray}
(The expression for the probability distribution of \eqref{*3} can
be written explicitly.)  Inequality \eqref{*2} follows from
\eqref{*1} due to the following property:
\begin{equation*}
\begin{aligned}
&\mathrm{Pr}\{\tau_{j,k}\leq x~|~\mathcal{A}_{i}\}\\
&:=\mathrm{Pr}\{\mbox{sum~of} \ i+1 \ \mbox{interarrival~times}\
\leq
x~|~\mbox{all~of~them~occur~during~a~service~time}\}\\
&\geq\mathrm{Pr}\{\mbox{sum~of} \ i+1 \ \mbox{interarrival~times}\
\leq x\}.
\end{aligned}
\end{equation*}
Inequality \eqref{*3} follows from \eqref{*2} due to the
elementary fact that for $i\leq j$
\begin{equation*}
\begin{aligned}
&\mathrm{Pr}\{\mbox{sum~of} \ i \ \mbox{interarrival~times}\ \leq
x\}\\
&\geq\mathrm{Pr}\{\mbox{sum~of} \ j \ \mbox{interarrival~times}\
\leq x\}.
\end{aligned}
\end{equation*}

Thus, $\tau_{j,k}$ is a proper random variable, and \eqref{r2.18}
is valid. Therefore, denoting the probability of \eqref{r2.18} by
$B_{j,k}(x)$, we see that under Conditions (A) and (D):
\begin{equation}\label{r2.18+}
\begin{aligned}
|B_{j,k}(x)-B(x)|&\leq\int_0^\infty
\underbrace{\sup_{x\geq0,y\geq0}\big|B_y(x)-B(x)\big|}_{<5\epsilon
\mbox{~due~to~Lemma \ref{l21}}}\mbox{d}\mathrm{Pr}\{\tau_{j,k}\leq
y~\}\\&<5\epsilon\int_0^\infty \mbox{d}\mathrm{Pr}\{\tau_{j,k}\leq
y~\}=5\epsilon,
\end{aligned}
\end{equation}
and \eqref{r2.18.1.new} follows.

Let $\kappa_{j,k}$ denote the number of inserted points within the
interval $[t_{j,k}, s_{j,k})$ \
($\sum_{k=1}^{f(j)}\kappa_{j,k}:=f(j+1)$). Then for $1\le j<n$
($n>1$)
\begin{equation}\label{r2.19}
\begin{aligned}
&\mathrm{Pr}\{\kappa_{j,k}=0\}=\int_0^\infty \mathrm{e}^{-\lambda
x}\mathrm{d}B_{j,k}(x)=\widehat{B}_{j,k}(\lambda),\\
&\mathrm{Pr}\{\kappa_{j,k}=m\}=[1-\mathrm{
Pr}\{\kappa_{j,k}=0\}][1-\widehat B(\lambda)]^{m-1} \widehat
B(\lambda),\quad m=1,2,...
\end{aligned}
\end{equation}
where $\widehat B(\lambda)=\int_0^\infty\mathrm{e}^{-\lambda
x}\mathrm{d}B(x)$. Notice, that according to \eqref{r2.18+} for
any $1\leq j<n$ ($n>0$)
\begin{equation}\label{r2.19+}
\begin{aligned}
|\widehat{B}_{i,j}(\lambda)-\widehat{B}(\lambda)|&=\left|\int_0^\infty\mathrm{e}^{-\lambda
x}\mathrm{d}B_{j,k}(x)-\int_0^\infty\mathrm{e}^{-\lambda
x}\mathrm{d}B(x)\right|\\
&\leq\int_0^\infty\lambda\mbox{e}^{-\lambda
x}\underbrace{|B_{i,j}(x)-B(x)|}_{<5\epsilon\mbox{~due~to~}\eqref{r2.18+}}\mbox{d}x<5\epsilon.
\end{aligned}
\end{equation}

The term $1-\mathrm{Pr}\{\kappa_{j,k}=0\}$ of \eqref{r2.19} is the
probability that during the residual service time of the tagged
customer there is at least one arrival (and therefore at least one
inserted point). The next term $[1-\widehat B(\lambda)]^{m-1}
\widehat B(\lambda)$ means that along with the first inserted
point associated with the tagged customer, there are $m-1$ other
inserted points.

For $j=n$ we have
\begin{equation}\label{r2.20}
\mathrm{Pr}\{\kappa_{n,k}=m\}=\int_0^\infty \mathrm{e}^{-\lambda
x}\frac{(\lambda x)^m}{m!}\mathrm{d}B_{n,k}(x),\quad m=0,1,...
\end{equation}

Let us now consider the time interval $[t_{0,1}, s_{0,1})$ and the
number of inserted points $\kappa_{0,1}$. Taking into account that
$B_{0,1}(x)=B(x)$ we have
\begin{equation}\label{r2.21}
\mathrm{Pr}\{\kappa_{0,1}=m\}=[1-\widehat B(\lambda)]^m \widehat
B(\lambda), \quad m=0,1,...,
\end{equation}
in the case when the loss queueing system is \textit{not}
$M/GI/1/0$, and
$$
\mathrm{Pr}\{\kappa_{0,1}=m\}=\int_0^\infty \mbox{e}^{-\lambda
x}\frac{(\lambda x)^m}{m!}\mathrm{d}B(x), \ m=0,1,...
$$
in the case of the $M/GI/1/0$ queueing system.

In the case when the queueing system is not $M/GI/1/0$,  for the
probability generating function of $\kappa_{0,1}$ from
(\ref{r2.21}) we have the representation
\begin{equation}\label{r2.23}
\sum_{m=0}^\infty z^m\mathrm{Pr}\{\kappa_{0,1}=m\}= \frac
{\widehat B(\lambda)}{1-z+z\widehat B(\lambda)}, \ \ |z|\le 1.
\end{equation}
Relation \eqref{r2.23} has been established in \cite{Abramov 2005}
as a particular case of a more general result.

Let us now prove the inequality \eqref{r2.8.1add}. For all $j<n$
and $k\geq1$ denote:
\begin{equation*}
\begin{aligned}
g_{j,k}(z)&=\sum_{m=0}^\infty z^m\mathrm{Pr}\{\kappa_{j,k}=m\},\ \
\
|z|\leq 1.\\
\end{aligned}
\end{equation*}

Comparing this expression with \eqref{r2.23} we have as follows.

Let us denote
$$
g^*(z)=\frac{\widehat B(\lambda)}{1-z+z\widehat B(\lambda)}.
$$
Then from \eqref{r2.19+} and \eqref{r2.23} we have
\begin{equation*}
\begin{aligned}
|g^*(z)-g_{j,k}(z)|&=\left|\sum_{m=0}^\infty z^m
\left[\mathrm{Pr}\{\kappa_{0,1}=m\}-\mathrm{Pr}\{\kappa_{j,k}=m\}\right]\right|\\
&=\left|\widehat B(\lambda)-\widehat
B_{j,k}(\lambda)+\sum_{m=1}^\infty z^m[1-\widehat
B(\lambda)]^{m-1}\widehat B^(\lambda)[\widehat{B}_{j,k}(\lambda)-\widehat B(\lambda)]\right|\\
&\leq\underbrace{|\widehat
B(\lambda)-\widehat{B}_{j,k}(\lambda)|}_{<5\epsilon \
\mbox{due~to~} \eqref{r2.19+}}+\underbrace{|\widehat
B(\lambda)-\widehat{B}_{j,k}(\lambda)|}_{<5\epsilon \
\mbox{due~to~} \eqref{r2.19+}}z\sum_{m=0}^\infty
z^m[1-\widehat B(\lambda)]^m\widehat B(\lambda)\\
&<5\epsilon[1+zg^*(z)].
\end{aligned}
\end{equation*}

Therefore,
\begin{equation*}
\begin{aligned}
g_{j,k}(z)&<5\epsilon[1+zg^*(z)]+g^*(z)\\
&=\frac {\widehat B(\lambda)}{1-z+z\widehat
B(\lambda)}+5\epsilon\cdot \frac {1-z+2z\widehat
B(\lambda)}{1-z+z\widehat B(\lambda)},
\end{aligned}
\end{equation*}
for all $j<n$ and $k\geq1$, and we arrive at
\begin{equation}\label{r2.24-1}
\begin{aligned}
g(z)=\sup_{j<n, k\geq1}g_{j,k}(z)&<\frac {\widehat
B(\lambda)}{1-z+z\widehat B(\lambda)}+5\epsilon\cdot \frac
{1-z+2z\widehat B(\lambda)}{1-z+z\widehat B(\lambda)},
\end{aligned}
\end{equation}
and \eqref{r2.8.1add} follows. Notice, that from the detailed
calculations of \eqref{r2.8}, inequality \eqref{r2.24-1} is
strong.

Now, in order to finish the proof of the theorem, let us first
establish some elementary properties of the class of probability
distribution $\{B_y(x)\}$. Recall that the explicit representation
for $B_y(x)$ follows from \eqref{r2.4-1}.

\begin{proper}\label{prop1}
The probability distribution $B_y(x)$ can be represented as
\begin{equation}\label{r2.100}
B_y(x)=r_yF_y(x)+(1-r_y)E_\mu(x),
\end{equation}
with parameter $r_y<1$ depending on $y$.
\end{proper}
The equivalent form of \eqref{r2.100} is
$$
\overline{B}_y(x)=r_y\overline{F}_y(x)+(1-r_y)\mathrm{e}^{-\mu x}.
$$

\begin{proof}
Indeed, we have:
$$
r_y\left(\frac{\overline{F}(x+y)}{\overline{F}(y)}\right)+(1-r_y)\mathrm{e}^{-\mu
x}=\frac{p\overline{F}(x+y)+(1-p)\mathrm{e}^{-\mu(x+y)}}{p\overline{F}(y)+(1-p)\mathrm{e}^{-\mu(y)}}.
$$
After some algebra we see that
$$
r_y=\frac{p\overline{F}(y)}{p\overline{F}(y)+(1-p)\mathrm{e}^{-\mu
y}},
$$
which clearly constitutes that $r_y<1$.
\end{proof}

We have the following property.
\begin{proper}
\label{prop2}$G\leq_{\mathscr{C}_{\lambda}} B_{y}$ for all
$y\geq0$.
\end{proper}

\begin{proof}
Indeed, since $G\leq_{\mathscr{C}_\lambda}F_{y^0}$,
$G\leq_{\mathscr{C}_\lambda}E_\mu$, then, because of
$F_{y^0}\leq_{\mathscr{C}_\lambda}F_{y}$, we also have
$G\leq_{\mathscr{C}_\lambda}F_{y}$ for all $y\geq0$. Therefore,
\begin{equation}\label{2.101}
r_{y}G+(1-r_{y})G\leq_{\mathscr{C}_{\lambda}}
r_{y}F_{y}+(1-r_{y})E_\mu.
\end{equation}
The left-hand side of \eqref{2.101} is $G(x)$, while the
right-hand side is $B_y(x)$. Therefore,
$G\leq_{\mathscr{C}_\lambda}B_y$.
\end{proof}

Property \ref{prop2} is crucial. Their consequence is the
following result $G\leq_{\mathscr{C}_\lambda}B_{n,k}$, where
$B_{n,k}(x)=\mathrm{Pr}\{\vartheta_{n,k}\leq x\}$, which follows
immediately from this Property \ref{prop2} by application of the
formula for the total probability. Indeed,
$$
B_{n,k}(x)=\mathrm{Pr}\{\vartheta_{n,k}\leq x\}=\int_0^\infty
B_y(x)\mathrm{d}\mathrm{Pr}\{\tau_{n,k}\leq y\},
$$
and consequently for $j=0,1,\ldots$,
\begin{equation*}
\begin{aligned}
\int_0^\infty\mathrm{e}^{-\lambda x}x^j
B_{n,k}(x)\mathrm{d}x&=\int_0^\infty\mathrm{e}^{-\lambda
x}x^j\left(\int_0^\infty
B_y(x)\mathrm{d}\mathrm{Pr}\{\tau_{n,k}\leq y\}\right)\mathrm{d}x\\
&\leq\int_0^\infty\mathrm{e}^{-\lambda
x}x^jG(x)\mathrm{d}x\int_0^\infty
\mathrm{d}\mathrm{Pr}\{\tau_{n,k}\leq y\}\\
&=\int_0^\infty\mathrm{e}^{-\lambda x}x^jG(x)\mathrm{d}x.
\end{aligned}
\end{equation*}

Hence, $G\leq_{\mathscr{C}_\lambda}B_{n,k}$.

Therefore, going back to (\ref{r2.20}) we conclude that, according
to Lemma A1 of Appendix A, there exists the minimal random
variable $\varsigma$ in the sense of the stochastic order
relation, i.e. $\varsigma\leq_{st}\kappa_{n,k}$ for all
$k=1,2,...$, corresponding to the minimal probability distribution
function $G(x)$ in the sense of the order relation
$\mathscr{C}_{\lambda}$. Specifically,
\begin{equation}\label{r2.43-added}
\mathrm{Pr}\{\varsigma=m\}=\int_0^\infty\mbox{e}^{-\lambda
x}\frac{(\lambda x)^m}{m!}\mathrm{d}G(x).
\end{equation}
Thus, taking the sequence of independent identically distributed
random variables $\varsigma_i$, $i=1,2,...$, all having the same
distribution as $\varsigma$, one can conclude
$$
L_n\geq_{st}\sum_{i=1}^{f(n)}\varsigma_i.
$$
Therefore, the statement of Theorem \ref{t24} follows due to
stochastic comparison.
\end{proof}

\begin{cor}\label{cor1}
In the case $p=\frac{1}{2}$ under Conditions (A) and (D) we have
$$
L_n\geq_{st}\sum_{i=1}^{Z_{n}}\varsigma_i.
$$
$Z_n$ denotes the number of offspring in the $n$th generation of
the Galton--Watson branching process with $Z_0=1$ and the
offspring generating function
$$
g(z)<\frac {\widehat B(\lambda)}{1-z+z\widehat
B(\lambda)}+2\epsilon\cdot \frac {1-z+2z\widehat
B(\lambda)}{1-z+z\widehat B(\lambda)}, \quad |z|\leq 1
$$
where
\begin{equation*}
\begin{aligned}
&\widehat B(\lambda)=p\widehat
F(\lambda)+(1-p)\frac{\mu}{\mu+\lambda},\\
&\widehat F(\lambda)=\int_0^\infty\mathrm{e}^{-\lambda
x}\mathrm{d}F(x).
\end{aligned}
\end{equation*}
 $\varsigma_i$, $i=1,2,\ldots$, is the sequence of nonnegative
integer-valued independent identically distributed random
variables all having the probability law
$$
\mathrm{Pr}\{\varsigma_i=m\}=\int_0^\infty\mathrm{e}^{-\lambda
x}\frac{(\lambda x)^m}{m!}\mathrm{d}G(x),
$$
and $G(x)$ is a probability distribution function that satisfies
the inequalities $G\leq_{\mathscr{C}_\lambda} F_{y^0}$ and
$G\leq_{\mathscr{C}_\lambda} E_\mu$, where
$F_{y^0}\leq_{\mathscr{C}_\lambda}F_{y}$ for all $y\ge 0$.
\end{cor}

\begin{proof} As we can see the only difference between the statements of
Theorem \ref{t24} and Corollary \ref{cor1} is in the coefficient
before $\epsilon$. Specifically, similarly to \eqref{r2.18+},
\begin{equation*}\label{r2.18++}
\begin{aligned}
|B_{j,k}(x)-B(x)|&\leq\int_0^\infty
\underbrace{\sup_{x\geq0,y\geq0}\big|B_y(x)-B(x)\big|}_{<2\epsilon\mbox{~due~to~}\eqref{r2.8_1-2}}
\mathrm{d}\mathrm{Pr}\{\tau_{j,k}\leq y~\}<2\epsilon,
\end{aligned}
\end{equation*}
and then similarly to \eqref{r2.19+}
\begin{equation*}
|\widehat{B}_{i,j}(\lambda)-\widehat{B}(\lambda)|<2\epsilon
\end{equation*}
as well, where
$$
B_{j,k}(x)=\mathrm{Pr}\{\vartheta_{j,k}\leq x\}.
$$
 The rest part of the proof of the corollary is similar to the
proof of Theorem \ref{t24}.
\end{proof}

\section{The number of losses under Condition (B)}

\noindent In this section it is assumed, additionally to Condition
(A) that the probability distribution function $F(x)$ belongs
either to the class NBU or to the class NWU. Having this
additional assumption, Condition (B) enables us to obtain stronger
inequalities than in the previous section under Condition (A).

\smallskip

The lemma below gives a stronger result than characterization
Lemma \ref{l22}. Specifically, we have the following lemma.

\begin{lem}\label{l31}
Under condition (B) we have
$$
\mathscr{K}(F,E_\mu)<\epsilon.
$$
\end{lem}

\begin{proof}
The proof of this lemma is given in Appendix B.
\end{proof}

\begin{lem}\label{l32}
Under Condition (B) we have
$$
\sup_{x\ge 0, y\ge 0}\big|B_y(x)-B(x)\big|<3\epsilon.
$$
\end{lem}

\begin{proof} We start the proof from representation (\ref{r2.4}).
Then, as in the proof of Lemma \ref{l21}, we  have (\ref{r2.5}).
However, taking into account that $F(x)$ belongs to one of the
aforementioned classes NBU and NWU, then instead of the earlier
inequality (\ref{r2.6}) given in the proof of Lemma \ref{l21} we
have the stronger inequality
\begin{equation}\label{r3.12}
|\overline{F}(x)-\mathrm{e}^{-\mu x}|<\epsilon,
\end{equation}
which in turn is the result of the application of Lemma \ref{l31}.
From (\ref{r2.4}), (\ref{r2.5}), (\ref{r3.12}) and the triangle
inequality, for all $x\ge 0$, $y\ge 0$ we obtain:
\begin{equation}\label{r3.13}
\begin{aligned}
&\big|p\overline{F}(x+y)-p^2\overline{F}(x) \overline{F}(y)
-p(1-p)\mathrm{e}^{-\mu
x}\overline{F}(y) \big|\\
&\leq p \overline{F}(y) \Big|\frac{ \overline{F}(x+y)}{
\overline{F}(y)}-\mathrm{e}^{-\mu
x}\Big|+p^2 \overline{F}(y) \big| \overline{F}(x)-\mathrm{e}^{-\mu x}\big|\\
&\leq p \overline{F}(y) \Big|\frac{ \overline{F}(x+y)}{
\overline{F}(y)}- \overline{F}(x)+
\overline{F}(x)-\mathrm{e}^{-\mu
x}\Big|\\&\ \ \ +p^2 \overline{F}(y) \big| \overline{F}(x)-\mathrm{e}^{-\mu x}\big|\\
&\leq p \overline{F}(y) \Big(\Big|\frac{ \overline{F}(x+y)}{
\overline{F}(y)}-
\overline{F}(x)\Big|+\big|\overline{F}(x)-\mathrm{e}^{-\mu
x}\big|\Big)\\
&\ \ \  +p^2 \overline{F}(y) \big| \overline{F}(x)-\mathrm{e}^{-\mu x}\big|\\
& < 2\epsilon p \overline{F}(y) +\epsilon p^2 \overline{F}(y)\\
&<3\epsilon p \overline{F}(y),
\end{aligned}
\end{equation}
\begin{equation}\label{r3.14}
\begin{aligned}
&\big|(1-p)\mathrm{e}^{-\mu(x+y)}-p(1-p)\mathrm{e}^{-\mu
y} \overline{F}(x) -(1-p)^2\mathrm{e}^{-\mu(x+y)}\big|\\
&\leq (1-p)\mathrm{e}^{-\mu y}\big|\mathrm{e}^{-\mu
x}-\overline{F}(x)\big|+(1-p)^2\mathrm{e}^{-\mu
y}\big|\mbox{e}^{-\mu
x}-\overline{F}(x)\big|\\
& < \epsilon(1-p)\mbox{e}^{-\mu y}+\epsilon(1-p)^2\mathrm{e}^{-\mu
y}\\
&<3\epsilon(1-p)\mathrm{e}^{-\mu y}.
\end{aligned}
\end{equation}
Therefore, from (\ref{r3.13}) and (\ref{r3.14}) for all $x\ge 0$,
$y\ge 0$ we have
$$
|B_y(x)-B(x)|<\frac{3\epsilon p \overline{F}(y) +3\epsilon
(1-p)\mathrm{e}^{-\mu y}}{ p\overline{F}(y)+(1-p)\mathrm{e}^{-\mu
y}}=3\epsilon.
$$
Lemma \ref{l32} is proved.
\end{proof}

\begin{rem}\label{rem2}
In the special case where $p=\frac{1}{2}$ we have:
\begin{equation}\label{r3.13_1}
\begin{aligned}
&\big|\textstyle{\frac{1}{2}} \overline{F}(x+y)
-\textstyle{\frac{1}{4}} \overline{F}(x) \overline{F}(y)
-\textstyle{\frac{1}{4}}\mathrm{e}^{-\mu
x} \overline{F}(y) \big|\\
&\leq
\textstyle{\frac{1}{4}} \overline{F}(y) \Big|\frac{ \overline{F}(x+y)}{ \overline{F}(y)}- \overline{F}(x)\Big|\\
&\ \ \ +\textstyle{\frac{1}{4}} \overline{F}(y) \Big|\frac{
\overline{F}(x+y)}{ \overline{F}(y)}-\mathrm{e}^{-\mu
x}\Big|\\
&\leq \textstyle{\frac{1}{2}} \overline{F}(y) \Big|\frac{
\overline{F}(x+y)}{ \overline{F}(y)}- \overline{F}(x)
\Big|\\
&\ \ \ +\textstyle{\frac{1}{4}} \overline{F}(y) \big|
\overline{F}(x)-\mathrm{e}^{-\mu x}\big|\\
&<\textstyle{\frac{3}{4}}\epsilon \overline{F}(y),
\end{aligned}
\end{equation}
\begin{equation}\label{r3.14_1}
\begin{aligned}
\big|\textstyle{\frac{1}{2}}\mathrm{e}^{-\mu(x+y)}-\textstyle{\frac{1}{4}}\mathrm{e}^{-\mu
y} \overline{F}(x)
-\textstyle{\frac{1}{4}}\mathrm{e}^{-\mu(x+y)}\big|
&=\textstyle{\frac{1}{4}}\mathrm{e}^{-\mu y}\big|\mathrm{e}^{-\mu
x}- \overline{F}(x)\big|\\
&<\textstyle{\frac{1}{4}}\epsilon\mathrm{e}^{-\mu y}\\
&<\textstyle{\frac{3}{4}}\epsilon\mathrm{e}^{-\mu y}.
\end{aligned}
\end{equation}
Therefore, from (\ref{r3.13_1}) and (\ref{r3.14_1}) we have:
$$
|B_y(x)-B(x)|<\frac{\frac{3}{4}\epsilon
 \overline{F}(y) +\frac{3}{4}\epsilon \mbox{e}^{-\mu
y}}{ \frac{1}{2}\overline{F}(y)+\frac{1}{2}\mbox{e}^{-\mu
y}}=\textstyle{\frac{3}{2}}\epsilon.
$$
\end{rem}

\medskip
The main result of this section is the following

\begin{thm}\label{t34}
Under Conditions (B) and (D) we have
$$
L_n\geq_{st}\sum_{i=1}^{Z_{n}}\varsigma_i.
$$
$Z_n$ denotes the number of offspring in the $n$th generation of
the Galton--Watson branching process with $Z_0=1$ and the
offspring generating function
$$
g(z)<\frac {\widehat B(\lambda)}{1-z+z\widehat
B(\lambda)}+3\epsilon\cdot \frac {1-z+2z\widehat
B(\lambda)}{1-z+z\widehat B(\lambda)}, \quad |z|\leq 1,
$$
where $\widehat B(\lambda)$ is as in Theorem \ref{t24}. The
sequence of independent identically distributed random variables
$\varsigma_i$, $i=1,2,\ldots$, is as in Theorem \ref{t24}.
\end{thm}

\begin{rem}
If Condition (D) is not satisfied, or the probability distribution
$F_{y^0}(x)$ is unknown and cannot be evaluated, then a
probability distribution function $G(x)$ should be chosen such
that $G\leq_{\mathscr{C}_\lambda} F_y$ for all $y\geq0$ and
$G\leq_{\mathscr{C}_\lambda} E_\mu$.
\end{rem}

\begin{proof} The proof of the theorem repeats the proof of corresponding
Theorem \ref{t24}. The only difference involves using the estimate
given by Lemma \ref{l32} rather than estimate given by Lemma
\ref{l21}.
\end{proof}

\begin{cor}\label{cor2}
In the case $p=\frac{1}{2}$ under Conditions (B) and (D) we have
$$
L_n\ge_{st}\sum_{i=1}^{Z_{n}}\varsigma_i.
$$
$Z_n$ denotes the number of offspring in the $n$th generation of
the Galton--Watson branching process with $Z_0=1$ and the
offspring generating function
$$
g(z)<\frac {\widehat B(\lambda)}{1-z+z\widehat
B(\lambda)}+\frac{3\epsilon}{2}\cdot \frac {1-z+2z\widehat
B(\lambda)}{1-z+z\widehat B(\lambda)}, \ \ |z|\le 1,
$$
where $\widehat B(\lambda)$ is as in Theorem \ref{t24}. The
sequence of independent identically distributed random variables
$\varsigma_i$, $i=1,2,\ldots$, is as in Theorem \ref{t24}.
\end{cor}

 The proof of Corollary \ref{cor2} is similar to the proof of Theorems
 \ref{t24} or \ref{t34}, or Corollary \ref{cor1}. The only
 difference is that the proof of this corollary uses the estimate obtained in Remark
 \ref{rem2} rather than the estimate of Lemma \ref{l32}.

\section{The number of losses under Condition (C)}
\noindent In this section we study the number of losses under
Condition (C). Applying Lemma \ref{l31} we have the following.

\begin{lem}\label{l41}
Under Condition (C)
$$
\mathscr{K}(B,E_\mu)<\epsilon.
$$
\end{lem}

The main result of the section is the following theorem.

\begin{thm}\label{t42}
Under Conditions (C) and (E) we have
\begin{equation}\label{rC1}
\sum_{i=1}^{X_{n}}\varsigma_i\leq_{st}L_n\leq_{st}\sum_{i=1}^{Y_{n}}\upsilon_i.
\end{equation}
$X_n$ denotes the number of offspring in the $n$th generation of
the Galton-Watson branching process with $X_0=1$ and the offspring
generating function
\begin{equation}\label{rC2}
g_X(z)<\frac {\widehat B(\lambda)}{1-z+z\widehat
B(\lambda)}+\epsilon\cdot \frac {1-z+2z\widehat
B(\lambda)}{1-z+z\widehat B(\lambda)}, \ \ |z|\le 1,
\end{equation}
and $Y_n$ denotes the number of offspring in the $n$th generation
of the Galton--Watson branching process with $Y_0=1$ and the
offspring generating function
$$
g_Y(z)= \frac {\widehat B(\lambda)}{1-z+z\widehat B(\lambda)}, \ \
|z|\leq 1.
$$
The sequence of independent identically distributed random
variables $\varsigma_i$, $i=1,2,\ldots$, is determined as follows:
\begin{equation}\label{eq1-4}
\mathrm{Pr}\{\varsigma_i=m\}=\int_0^\infty\mathrm{e}^{-\lambda
x}\frac{(\lambda x)^m}{m!}\mathrm{d}B_{y^0}(x), \quad
m=0,1,\ldots,
\end{equation}
where $B_{y^0}(x)$ is the minimal probability distribution
function in the sense that $B_{y^0}\leq_{\mathscr{C}_\lambda}B_y$
for all $y\ge 0$. The sequence of independent identically
distributed random variables $\upsilon_i$, $i=1,2,\ldots$, in turn
is determined as follows:
\begin{equation}\label{eq2-4}
\mathrm{Pr}\{\upsilon_i=m\}=\int_0^\infty\mathrm{e}^{-\lambda
x}\frac{(\lambda x)^m}{m!}\mathrm{d}B(x), \quad m=0,1,\ldots
\end{equation}
\end{thm}

\begin{rem}\label{rem-5-}
If the probability distribution $B_{y^0}(x)$ does not exist, then
we only have the one-side stochastic inequality $
L_n\leq_{st}\sum_{i=1}^{Y_{n}}\upsilon_i $.
\end{rem}

\begin{proof} Notice that under the assumption that
$B(x)$ belongs to the class NBU, we have the inequality
$\sum_{i=1}^{Y_{n}}\upsilon_i\geq_{st}L_n$, where $\upsilon_i$,
$i=1,2,\ldots$, is a sequence of independent identically
distributed random variables, defined by the probability law $
\mathrm{Pr}\{\upsilon_i=m\}=\int_0^\infty\mathrm{e}^{-\lambda
x}\frac{(\lambda x)^m}{m!}\mathrm{d}B(x) $.

 Indeed, under this
assumption from (\ref{r2.19}) and (\ref{r2.21}) we obtain that
$\kappa_{j,k}\leq_{st}\kappa_{0,1}$ for all $0\leq j<n$ ($n>0$)
and $k\ge 1$ (for details see Abramov \cite{Abramov 2005}), and
then the desired inequality
$\sum_{i=1}^{Y_{n}}\upsilon_i\ge_{st}L_n$ follows.

In order to prove the second inequality
$L_n\ge_{st}\sum_{i=1}^{X_{n}}\varsigma_i$ let us note the
following. Assumption (\ref{r30}) means that for all $0\leq j<n$
($n>0$), $k\geq 1$ we have
$\mathscr{K}(\vartheta_{j,k},\chi)<\epsilon$. Then using the same
arguments as in the proof of Theorem \ref{t24} we arrive at the
left-side of inequality \eqref{rC1} with the offspring generating
function defined by inequality \eqref{rC2}. The remaining part of
the proof of the theorem is similar to the corresponding part of
the proof of Theorem \ref{t24}.
\end{proof}

\section*{Appendix A: The order relation $\mathscr{C}_\lambda$}

Let $\xi_1$, $\xi_2$ be two nonnegative random variables, and let
$\Xi_1(x)$ and $\Xi_2(x)$ be their probability distribution
functions respectively.

\medskip
\noindent \textbf{Definition A1.} The random variable $\xi_1$ is
said to be less than the random variable $\xi_2$ in the sense of
the relation $\mathscr{C}_\lambda$ (notation:
$\xi_1\le_{\mathscr{C}_{\lambda}}\xi_2$ or
$\Xi_1\le_{\mathscr{C}_{\lambda}}\Xi_2$) if for fixed parameter
$\lambda>0$ and all $i=0,1,...$ there is the inequality
$$
\int_0^\infty\mbox{e}^{-\lambda
x}x^i\Xi_1(x)\mbox{d}x\ge\int_0^\infty\mbox{e}^{-\lambda
x}x^i\Xi_2(x)\mbox{d}x.
$$

Relation $\mathscr{C}_\lambda$ is a partial order relation. It
follows from the following lemma.

\medskip
\noindent \textbf{Lemma A1}. {\it Let
$\xi_1\leq_{\mathscr{C}_{\lambda}}\xi_2$. Then for the random
variables $\theta_1$ and $\theta_2$ given by the probability
laws:}
$$
\mathrm{Pr}\{\theta_k=m\}=\int_0^\infty\mbox{e}^{-\lambda
x}\frac{(\lambda x)^m}{m!}\mathrm{d}\Xi_k(x), \ k=1,2; \
m=0,1,...,
$$
{\it we have $\theta_1\leq_{st}\theta_2$. Converse, if
$\theta_1\leq_{st}\theta_2$, then
$\xi_1\leq_{\mathscr{C}_{\lambda}}\xi_2$.}

\begin{proof} Indeed, by partial integration we have
$$
\int_0^\infty\mbox{e}^{-\lambda
x}\mathrm{d}\Xi_k(x)=\lambda\int_0^\infty\mbox{e}^{-\lambda x}
\Xi_k(x)\mathrm{d}x \quad (k=1,2),
$$
and
\begin{equation*}
\begin{aligned}
&\int_0^\infty\mathrm{e}^{-\lambda x}\frac{(\lambda
x)^m}{m!}\mbox{d}\Xi_k(x)\\
&= \lambda\Big(\int_0^\infty\mathrm{e}^{-\lambda x}\frac{(\lambda
x)^m}{m!}\Xi_k(x)\mathrm{d}x-\int_0^\infty\mbox{e}^{-\lambda x}
\frac{(\lambda x)^{m-1}}{(m-1)!}\Xi_k(x)\mathrm{d}x\Big),
\end{aligned}
\end{equation*}
$$
(m=0,1,...; \ k=1,2).
$$
Therefore, for all $m=0,1,...,$ and $k=1,2$, we obtain
$$
\lambda\int_0^\infty\mbox{e}^{-\lambda x}\frac{(\lambda
x)^m}{m!}\Xi_k(x)\mbox{d}x=\sum_{i=0}^m\int_0^\infty\mathrm{e}^{-\lambda
x}\frac{(\lambda x)^i}{i!}\mathrm{d}\Xi_k(x),
$$
and the result follows.
\end{proof}

\section*{Appendix B: The proof of Lemmas \ref{l22} and \ref{l31}}
The only difference between the proofs of these lemmas is that
there is an additional condition (B) in Lemma \ref{l31}.
Therefore, we provide the general proof of both Lemmas \ref{l22}
and \ref{l31} with an appropriate specification of the case
studies.

We have
$$
\overline{F}(x+y)=\overline{F}(y) \overline{F}(x)+h(x,y),
$$
where $h(x,y)=F(x)-\mathrm{Pr}\{\zeta\leq x+y~|~\zeta>y\}$. Then,
$$ \int_0^\infty
\overline{F}(x+y)\mathrm{d}y=\frac{1}{\mu}\overline{F}(x)+\frac{1}{\mu}\Theta(x),\leqno(B.1)
$$
where
$$
\Theta(x)=\mu\int_0^\infty \overline{F}(y)
h(x,y)\mathrm{d}y.\leqno(B.2)
$$
Denoting
$$
\Psi(x)=\int_0^\infty \overline{F}(x+y)\mathrm{d}y=\int_x^\infty
\overline{F}(y) \mathrm{d}y,
$$
from (B.1) we obtain the following differential equation
$$
\Psi'(x)+\mu\Psi(x)-\Theta(x)=0.\leqno(B.3)
$$
The solution of differential equation (B.3), satisfying the
initial condition $\Psi(0)=1/\mu$, is
$$
\Psi(x)=\frac{1}{\mu}\mbox{e}^{-\mu
x}+\int_0^x\mathrm{e}^{-\mu(x-t)}\Theta(t)\mathrm{d}t.
$$
Hence,
$$
\overline{F}(x)=\mbox{e}^{-\mu
x}+\mu\int_0^x\mbox{e}^{-\mu(x-t)}\Theta(t)\mbox{d}t-\Theta(x).\leqno(B.4)
$$
Let us study some properties of (B.4). Note first that
$|\Theta(x)|<\epsilon$. Indeed, according to assumption of this
lemma $\sup_{x\ge 0, y\ge 0}|h(x,y)|<\epsilon$, and therefore, as
it follows from (B.2),
$$
|\Theta(x)|\leq\Big(\sup_{x\ge 0, y\geq
0}\big|h(x,y)\big|\Big)\mu\int_0^\infty
\overline{F}(y)\mathrm{d}y<\epsilon.
$$
Therefore,
$$
\Big|\mu\int_0^x\mathrm{e}^{-\mu(x-t)}\Theta(t)\mathrm{d}t-\Theta(x)\Big|\leq
\Big|\mu\int_0^x\mathrm{e}^{-\mu(x-t)}\Theta(t)\mathrm{d}t\Big|+|\Theta(x)|
$$
$$
\leq
2\sup_{0\leq t\leq x}|\Theta(t)|<2\epsilon.\leqno(B.5)
$$
Substituting (B.5) for (B.3), for all $x\ge 0$ we obtain
$$
|\overline{F}(x)-\mathrm{e}^{-\mu x}|<2\epsilon,
$$
and the statement of Lemma \ref{l22} follows.

Under the assumption that $F(x)$ belongs either to the class NBU
or to the class NWU, the function $\Theta(x)$ is either
nonnegative or nonpositive. This is because $h(x,y)$ contains the
term belonging to one of these classes, and therefore is either
nonnegative or nonpositive. Therefore, according to (B.2)
$$
\Big|\mu\int_0^x\mbox{e}^{-\mu(x-t)}\Theta(t)\mbox{d}t
-\Theta(x)\Big|\leq|\Theta(x)|<\epsilon.\leqno(B.6)
$$
Substituting (B.6) for (B.3), for all $x\ge 0$ we obtain
$$
|\overline{F}(x)-\mbox{e}^{-\mu x}|<\epsilon,
$$
and Lemma \ref{l31} is therefore proved as well.

\section*{Appendix C: The proof of Lemma \ref{tex-lem}}
Let $\mathscr{T}:=\{t_n, n\geq0\}$ denote the time instants at
which customers arrive to the system (a customer can arrive
without joining queue), where $t_0=0$, and let $\{N(C),
C\in\mathscr{B}(\mathbb{R})\}$ denote the Poisson process that
consists of the collection of points $\mathscr{T}$.

Next, let $\{{T}_n, n\geq1\}$ denote the time instants at which a
new service begins in the system. At these time instants either

\smallskip
(i) a departure occurs and leaves at least one customer behind the
system, or

\smallskip
(ii) an arrival occurs and finds an empty system.

\smallskip

Finally, let $Q(t)$ denote the number of customers in the system
at time $t$.

Let us find $\mathrm{Pr}\{\vartheta_{j,k}\leq x|\tau_{j,k}=y\}$.
Suppose $\xi_m=t_{N(T_m)+j-Q(T_m)}-T_m$. Clearly, $\xi_m$ is the
remaining amount of time, after the beginning of the $m$th
service, until a new customer arrives, plus the next $j-Q(T_m)-1$
interarrival times. The instant $T_m$ is chosen such that
$s_{j,k-1}\leq T_m<t_{j,k}$, where we set $s_{j,0}:=0$. (The
definition of the time instants $s_{j,k}$ and $t_{j,k}$ is given
in the proof of Theorem \ref{t24}.) Denoting $\eta=\inf\{n\geq1:
\chi_m>\xi_m\}$, where $\chi_m$ denote the length of the $m$th
service time, one can see that
$$
\mathrm{Pr}\{\vartheta_{j,k}\leq
x|\tau_{j,k}=y\}=\mathrm{Pr}\{\chi_\eta-\xi_\eta\leq
x|\xi_\eta=y\},
$$
and
\begin{equation*}
\begin{aligned}
&\mathrm{Pr}\{\chi_\eta-\xi_\eta\leq x|\xi_\eta=y\}\\
&=\sum_{n=1}^\infty\mathrm{Pr}\{\chi_\eta-\xi_\eta\leq x|\eta=n,
\xi_n=y\}\mathrm{Pr}\{\eta=n|\xi_\eta=y\}\\
&=\sum_{n=1}^\infty\mathrm{Pr}\{\chi_n-\xi_n\leq x|\chi_m<\xi_m,
1\leq m\leq n-1, \chi_n>y, \xi_n=y\}\mathrm{Pr}\{\eta=n|\xi_\eta=y\}\\
&=\sum_{n=1}^\infty\frac{\mathrm{Pr}\{\chi_n\leq x+y,
\chi_n>y|\chi_m<\xi_m, 1\leq m\leq n-1, \xi_n=y\}}
{\mathrm{Pr}\{\chi_n>y|\chi_m<\xi_m, 1\leq m\leq n-1,
\xi_n=y\}}\mathrm{Pr}\{\eta=n|\xi_\eta=y\}.
\end{aligned}
\end{equation*}
Notice, that each $\xi_n$ is a function of the Poisson arrival
process and the random variables $\chi_1$, $\chi_2$, \ldots,
$\chi_{n-1}$, and therefore the random variable $\chi_n$ is
independent of the vector ($\xi_1$, $\chi_1$, $\xi_2$,
$\chi_2$,\ldots, $\xi_{n-1}$, $\chi_{n-1}$, $\xi_n$), and so
$$
\frac{\mathrm{Pr}\{\chi_n\leq x+y, \chi_n>y|\chi_m<\xi_m, 1\leq
m\leq n-1, \xi_n=y\}} {\mathrm{Pr}\{\chi_n>y|\chi_m<\xi_m, 1\leq
m\leq n-1, \xi_n=y\}}=B_y(x).
$$
Hence,
$$
\mathrm{Pr}\{\chi_\eta-\xi_\eta\leq x|\xi_\eta=y\}=B_y(x),
$$
or
$$
\mathrm{Pr}\{\vartheta_{j,k}\leq x|\tau_{j,k}=y\}=B_y(x).
$$
and the statement of Lemma \ref{tex-lem} is proved.

\section*{Acknowledgements}
The author thanks the anonymous referees for their comments. Owing
their comments and advice the paper has been changed
significantly. Author thanks Prof. Kais Hamza (Monash University)
and Prof. Robert Liptser (Tel Aviv University) for useful
conversation and Prof. Daryl Daley (Australian National
University) for useful editorial comments and conversation. An
on-line question of Taleb Samira (USTHB, Algeria) led to a
correction of a flaw. The paper is dedicated to V.M. Zolotarev,
V.M. Kruglov and V.V. Kalashnikov. These Academics helped the
author essentially, when he lived in the former USSR. The author
also acknowledges with thanks the support of the Australian
Research Council (grant \#DP0771338).


\end{document}